\documentclass[12pt,leqno,a4paper,dvips,twoside]{article}
\usepackage{ngerman}
\usepackage{amsfonts}
\usepackage{graphicx}
\usepackage{amsmath,amsfonts,amssymb}
\usepackage{makeidx}
\usepackage{shortvrb}
\usepackage{empheq}
\usepackage{extarrows}

\usepackage{graphicx}
\usepackage{amsfonts}
\usepackage{amssymb}
\usepackage{amsmath}
\usepackage[T1]{fontenc}

\newtheorem{theo}{Theorem}

\begin{document}

\begin{center}
{\large {\bf  S\"atze vom {\sc Bohman-Korovkin}-Typ f\"ur lokalkonvexe Vektorverb\"ande}}
\vspace{0.5cm}

{\large Heiner Gonska }
\end{center}
\vspace{0.5cm}


\bigskip \bigskip




\section{Einleitung}

Wir betrachten das Schema
$$
C(X) \xlongrightarrow[\text{P}]{\text{$T_n$}} F .
$$
Hierin bezeichnet $C(X)$ den Raum der stetigen reellwertigen Funktionen auf einem kompakten Raum $X$ und $F$ einen lokalkompakten Vektorverband. $(T_n)_{n \in \mathbb{N}}$ bezeichne eine Folge positiver linearer Abbildungen und $P$ einen Verbandshomomorphismus. In dieser Situation gilt der folgende Satz, der eine Verallgemeinerung eines Ergebnisses von Berens und Lorentz \cite{BeLo} darstellt.

\bigskip \bigskip

\noindent {\bf Satz A}

\medskip

Sei $F$ ein lokalkonvexer Vektorverband und $P : C(X) \to F$ ein Verbandshomomorphimus. Ist $S$ eine Teilmenge von $C(X)$, die eine strikt positive Funktion $g^*$ enth\"alt und gilt f\"ur eine Folge positiver Abbildungen $T_n : C(X) \to F$, $n \in \mathbb{N}$,
$$
T_ng \to Pg \;\; \mbox{in}\;\; F \;\; \mbox{f\"ur alle}\;\; g \in G = \mbox{ lin } S ,
$$
so folgt
$$
T_ng \to Pg \;\; \mbox{in}\;\; F \;\; \mbox{f\"ur alle}\;\; f \in \hat{G}_{supp P} .
$$
Dabei ist
$$
\begin{array}{lcl}
\hat{G}_{supp P} & = & \left\{ f \in C(X) : \mu (f) = f(x) \;\; \mbox{f\"ur alle} \;\; x \in supp P \right. \;\; \\
 & &  \mbox{und alle positiven Linearformen} \;\; \mu \;\; \mbox{mit}\;\; \mu (g) = g (x) \;\; \\
 & &  \left.\mbox{f\"ur alle}\;\; g \in G\right\} \\
\end{array}
$$
der Fortsetzungsraum bzgl. $supp P$ und $G$.
F\"ur einen Beweis siehe \cite{Gon75}, Theorem 3.1. Hieraus ergibt sich unmittelbar

\bigskip

\noindent {\bf Satz B} \\

\noindent Es sei $F$ ein lokalkonvexer  Vektorverband und $P : C(X) \to F$ ein Verbandshomomorphismus. Ist $S$ eine Teilmenge von $C(X)$, die eine strikt positive Funktion $g^*$ enth\"alt und $G = \mbox{ lin } S$, so gilt
$$
\hat{G}_{supp P} \subset \rho (S,F,P) .
$$
Hierzu bezeichnet $\rho (S,F,P)$ den Schatten von $S$ bzgl. $L^+ (C(X),F)$ - der Menge aller positiven linearen Abbildungen von $C(X)$ nach $F$ - und $P$, d.h.
$$
\begin{array}{lcl}
\rho (S,F,P) & = & \{ f\in C(X) : \;\; \mbox{Ist} \;\; (T_n)_{n \in \mathbb{N}} \subset L^+ (C(X)F) \;\; \mbox{mit} \;\; T_ng \to Pg \\
& & \mbox{f\"ur alle} \;\; g \in S, \;\; \mbox{so folgt} \;\; T_nf \to Pf\}.
\end{array}
$$

\section{Ergebnisse}

Wir beginnen mit einer Verallgemeinerung eines Satzes von {\sc Scheffold} \cite{scheffold}.

\begin{theo}
Es sei $E$ ein topologischer Vektorraum und ein Vektorverband und $F$ ein lokalkonvexer Vektorverband, $X$ sei eine kompakte Menge und $T: C(X) \to E$ ein Verbandshomomorphismus. 
Es sei $T_n : E \to F$ eine gleichstetige Folge positiver linearer Abbildungen und $P : E \to F$ ein stetiger Verbandshomomorphismus. Es sei $S$ eine Teilmenge von $C(X)$, die eine strikt positive Funktion $g*$ enth\"alt und $G = lin S$. 

Es gelte
$$
T_n (Tg) \to P(Tg) \;\;Ê\mbox{in} \;\; F \;\; \mbox{f\"ur alle} \;\; g \in S.
$$
Dann folgt
$$
T_n(\overline{f}) \to P(\overline{f}) \;\; \mbox{in} \;\; F \;\; \mbox{f\"ur alle} \;\; \overline{f} \in \overline{T(\hat{G}_{supp P\circ T})}^E.
$$

\end{theo}

\bigskip

\noindent{\bf Beweis:} $(T_n \circ T)_{n \in \mathbb{N}}$ ist eine Folge positiver linearer Abbildungen mit
$$
T_n \circ T : C(X) \to F.
$$
$P \circ T : C(X) \to F$ ist ein Verbandshomomorphismus. 

Nach Satz B gilt
$$
\hat{G}_{supp P \circ T} \subset \rho (S,F,P \circ T) .
$$
Aus
$$
T_n \circ T(g) \to P \circ T(g) \;\; \mbox{in} \;\; F \;\; \mbox{f\"ur alle} \;\; g \in S ,
$$
folgt also
$$
T_n \circ T(h) \to P \circ T(h) \;\; \mbox{in} \;\; F \;\; \mbox{f\"ur alle} \;\; h \in \hat{G}_{supp P \circ T},
$$
also
$$
T_n (f) \to P(f) \;\; \mbox{in} \;\; F \;\; \mbox{f\"ur alle} \;\; f \in T(\hat{G}_{supp P \circ T}).
$$
Sei nun $\overline{f}  \in \overline{T(\hat{G}_{supp P \circ T})}^E$ und $U$ eine beliebige Nullumgebung in $F$. Dann existiert eine Nullumgebung $W$ in $F$ mit
$$
W + W + W \subset U .
$$
Wegen der Gleichstetigkeit der Folge $(T_n)_{n \in \mathbb{N}}$ existiert eine Umgebung $V$ in $E$ mit $T_n (V) \subset W$ f\"ur alle $n \in \mathbb{N}$ und wegen der Stetigkeit von $P$ eine Umgebung $V'$ mit $P(V') \subset W$. Es sei nun $h \in \hat{G}_{supp P \circ T}$ so gew\"ahlt, dass $\overline{f} - T(h) \in V \cap V'$ ist. F\"ur $n \ge N_0$ ist dann
$$
T_n \circ T(h) - P \circ T(h) \in W,
$$
und insgesamt ergibt sich
{\footnotesize
$$
\begin{array}{cccccccccccccccc}
T_n\overline{f} & - & P \overline{f} & = & T_n \overline{f} & - & T_n \circ T(h) & + & T_n \circ T(h) & - & P \circ T(h) & + & P \circ T(h) & - & P\overline{f}  & \\[2mm]
                 
                &   &                & \in &                & W &                & + &                & W &              & + &              & W &                \subset U 
\end{array}
$$
}
f\"ur alle $n \ge N_0$. \hfill $\square$

\bigskip \bigskip

\noindent{\bf Bemerkung:} {\sc Scheffold} \cite{scheffold} hat Theorem 1 f\"ur einen lokalkonvexen Vektorverband $E$ und einen injektiven Verbandshomomorphimus $T : C(X) \to E$ 
mit v\"ollig anderen Mitteln bewiesen. Statt von Konvergenz auf dem Abschluss des Bildes eines Fortsetzungsraumes zu sprechen, verwendet {\sc Scheffold} die Terminologie
des relativen Choquetrandes.

\bigskip

Nach Theorem 1 sind nat\"urlich solche Vektorverb\"ande und topologische Vektorr\"aume $E$ von Interesse, in die ein Raum $C(X)$ verm\"oge einer nat\"urlichen Inklusion eingebettet ist und f\"ur die gilt $\overline{i(C(X))}^E = E$. Eine Klasse solcher R\"aume bilden sogenannte Banachsche Funktionenr\"aume. Dazu die

\bigskip \bigskip

\noindent {\bf Definition 2} {\sc (M\"uller)} \cite{mueller}) \\
Es sei $K \subset \mathbb{R}$ kompakt und $M(K)$ der Vektorraum der auf $K$ definierten reellwertigen (Lebesgue-)me\ss{}baren Funktionen modulo des zugeh\"origen Nullraumes $N$. Ein Banach-Raum $(B (K) , ||\circ||_B)$ bestehend aus Elementen von $M(K)$ hei\ss{}t ein {\it Banachscher Funktionenraum} genau dann, wenn seine Norm den folgenden Bedingungen gen\"ugt:
\begin{itemize}
\item[(N1)] Ist $g \in M(K)$ und $f \in B(K)$ mit $|g| \le |f|$, so folgt: $g \in B(K)$ und $||g||_B \le ||f||_B$.
\item[(N2)] Ist $(f_n)_{n \in \mathbb{N}} $ eine Folge in $B(K)$ und $0 \le f_n \nearrow f$ mit $f \in B(K)$, so folgt: $||f_n||_B \to ||f||_B$. 
\item[(N3)] $||f||_B$ ist umordnungsinvariant f\"ur alle $f \in B(K)$, d.h.: Ist $f = f'$ $\lambda$-fast \"uberall, so ist $f' \in B(K)$ und $||f||_B = ||f'||_B$.
\end{itemize}

\bigskip

\noindent {\bf Beispiel} {\sc (M\"uller)} \cite{mueller}) \\
Der Raum $L^p(K), \;\; 1 \le p < \infty$ ist ein Banachscher Funktionenraum.

\bigskip \bigskip

\noindent F\"ur Banachsche Funktionenr\"aume, die den Raum $C(K)$ enthalten, gilt der

\bigskip \bigskip

\noindent {\bf Satz 3} {\sc (M\"uller)} \cite{mueller}) \\
Ist $K \subset \mathbb{R}$ und $B(K)$ ein Banachscher Funktionenraum, der den Raum $C(K)$ enth\"alt, so ist $C(K)$ dicht in $B(K)$.

\bigskip

\noindent Um Theorem 1 anwenden zu k\"onnen, ben\"otigen wir noch

\bigskip

\noindent {\bf Lemma 4} {\it Ein Banachscher Funktionenraum $B(K)$ ist ein Banachverband.}

\noindent {\bf Beweis:} Es seien $f $ und $g \in B(K)$.
Dann ist $\sup (f,g) = \frac{f+g}{2} + \frac{|f-g|}{2}$. 
Wegen (N1) sind $|f|, |g| \in B(K)$ und wegen $|f-g| \le |f| + |g|$ ist auch $|f-g| \in B(K)$ und damit auch $\sup (f,g)$. Da $\inf (f,g) = \frac{f+g}{2} - \frac{|f-g|}{2}$ gilt, folgt, dass auch $\inf (f,g) \in B(K)$ ist. $B(K)$ ist also ein Vektorverband.

\bigskip

Eine Nullumgebungsbasis in $B(K)$ bilden die Normkugeln $K(0,\epsilon) = \{ f \in B(K) : ||f||_B \le \epsilon \}$. 
Ist nun $g \in B(K)$ mit $|g| \le |f|$ und $||f||_B \le \epsilon$, so folgt mit (N1), dass $||g||_B \le ||f||_B \le \epsilon$, also $g \in K(0,\epsilon)$ ist. Dies bedeutet aber, dass $K(0,\epsilon) $ solide ist, und der Banachraum und Vektorverband $B(K)$ besitzt eine Nullumgebungsbasis aus soliden Mengen. $B(K)$ ist also ein Banachverband. \hfill $\square$

\bigskip \bigskip

\noindent Wir beweisen nun eine Verallgemeinerung eines Satzes von M\"uller {\sc (M\"uller)} \cite{mueller}).

\bigskip \bigskip

\noindent {\bf Satz 5} \\
Es sei $K \subset \mathbb{R}$ kompakt und $B(K)$ ein Banachscher Funktionenraum, der den Raum $C(K)$ enth\"alt. Es sei $S$ eine Teilmenge von $C(K)$, die eine strikt positive Funktion $g^*$ enth\"alt und $G = \mbox{lin} S$. Es gelte $\hat{G}_K = C(K)$. 
Es sei $T_n : B(K) \to B(K)$, $n \in \mathbb{N}$, eine Folge positiver linearer Operatoren mit
$$
\sup_{n \in \mathbb{N}} ||T_n|| < + \infty
$$
und
$$
\lim_{n \to \infty} ||T_ng - g||_B = 0 \;\; \mbox{f\"ur alle}\;\; g \in S .
$$
Dann gilt:
$$
\lim_{n \to \infty} ||T_nf - f||_B = 0 \;\; \mbox{f\"ur alle} f \in B(K) .
$$

\noindent {\bf Beweis:} Die Injektion $i : C(K) \ni f \mapsto f \in B(K)$ ist ein Verbandshomomorphismus. Wegen $\sup_{n \in \mathbb{N}} ||T_n|| < + \infty$ ist die Folge $(T_n)_{n\in \mathbb{N}}$ gleichstetig. Die Identit\"at $I : B(K) \ni f \mapsto f\in B(K)$ ist ein stetiger Verbandshomomorphismus.

\bigskip

\noindent Mit Theorem 1 folgt also
$$
T_n f \to If \;\; \mbox{in} \;\; B(K) \;\;Ê\mbox{f\"ur alle} \;\; f \in \overline{i({\hat{G}}_{supp \; I \circ i})}^{B(K)}.
$$
Nun ist
$$
i(\hat{G}_{supp \; I \circ i}) \supset i(\hat{G}_K) = i(C(K)) = C(K) ,
$$
also 
$$
\overline{i(\hat{G}_{supp \; I \circ i})}^{B(K)} = \overline{C(K)}^{B(K)} = B(K). 
$$
\mbox{} \hfill $\square$

\bigskip

Mit Satz 5 ist es nun etwa m\"oglich, den folgenden Approximationsprozess im Banachschen Funktionenraum $L^1 [a,b]$, $[a,b] \subset \mathbb{R}$, auf Konvergenz gegen die Identit\"at zu testen. 

\bigskip \bigskip

\noindent {\bf Beispiel 6} \\
F\"ur $n \in \mathbb{N}$ und $f \in L^1 [a,b]$ sei
$$
T_n (f,x) := \frac{1}{b-a} \circ \sqrt{\frac{n}{\pi}} \circ \int_a^b \left( 1-\left( \frac{t-x}{b-a} \right)^2 \right)^n \circ f(t) dt.
$$
Darin hei\ss{}t
$$
K_n (t,x) := \frac{1}{b-a} \circ \sqrt{\frac{n}{\pi}} \circ \left( 1 - \left( \frac{t-x}{b-a} \right)^2 \right)^n
$$
 Landau-Stieltjes-Kern.\\
$T_n$ ist f\"ur alle $n \in \mathbb{N}$ linear und positiv, und f\"ur alle $f \in L^1 [a,b]$ ist $T_n f \in \Pi_{2n} \subset L^1 [a,b]$. Wie {\sc M\"uller} \cite{mueller} zeigt, sind die Normen der $T_n$ gleichm\"a\ss{}ig beschr\"ankt und f\"ur $0 \le i \le 2$ gilt $T_n \pi_i \to \pi_i$ in $L^1 [a,b]$. F\"ur $G = lin \{\pi_0,\pi_1,\pi_2\}$ gilt  $\hat{G}_K = C(K)$, so dass mit Satz 5 nun
$$
T_n f \to f \;\; \mbox{in}\;\; L^1 [a,b] \;\; \mbox{f\"ur alle}\;\; f \in L^1 [a,b]
$$
folgt. \hfill $\square$

\bigskip

Wir betrachten nun die spezielle Klasse lokalkonvexer Vektorverb\"ande $E$, die folgenvollst\"andig sind und im positiven Kegel einen quasiinneren Punkt besitzen. Es gilt folgendes

\bigskip

\noindent {\bf Theorem 7} (vgl. Scheffold \cite{scheffold}) {\it Es sei $E$ ein lokalkonvexer, folgenvollst\"andiger Vektorverband, $u \in E$ ein quasiinnerer Punkt des positiven Kegels $E^+$ von $E$ und $\{u_i , i \in I\}$ ein Erzeugendensystem von $E$ mit
$$
\{u_i , i \in I\} \subset E_u = \bigcup_{n =1}^\infty \{ x \in E : |x| \le n \circ u\}.
$$
Dann gilt:
\begin{itemize}
\item[(i)] Es gibt einen kompakten Raum $X$ und einen surjektiven Verbandsisomorphismus $T : C(X) \to E_u$.
\item[(ii)] Setzt man $u_i^{(2)} := T((T^{-1} (u_i))^2)$, so gilt: Ist $F$ ein beliebiger lokalkonvexer Vektorverband, $P : E \to F$ ein stetiger Verbandshomomorphismus und $T_n : E \to F$, $n \in \mathbb{N}$, eine gleichstetige Folge positiver linearer Abbildungen mit
$$
T_n v \to Pv \;\; \mbox{in} \;\; F \;\; \mbox{f\"ur alle} \;\; v \in \{u\} \cup \{u_i,u_i^{(2)}; i \in I\} ,
$$
so folgt
$$
T_nf \to Pf \;\; \mbox{in} \;\; F \;\; \mbox{f\"ur alle} \;\; f \in E .
$$
\item[(iii)] Wird $E$ von $k$ Elementen erzeugt, so kann die Kovnergenz einer Folge gleichstetiger positiver linearer Abbildungen gegen einen stetigen Verbandshomomorphismus auf einer h\"ochstens $(2k + 1)$-elementigen Menge getestet werden.
\end{itemize}
}

\bigskip

\noindent {\bf Beweis:} Nach Definition des quasiinneren Punktes $u$ ist das Vektorverbandsideal $E_u = \bigcup_{n=1}^\infty \{x \in E : |x| \le n \circ u\}$ dicht in $E$.
Damit existieren also ein kompakter Hausdorffraum $X$ und ein surjektiver Verbandsisomorphismus $T : C(X) \to E_u$ mit $T1_X = u$. Also gilt die Aussage (i).

\bigskip

Wir betrachten nun
$$
Q := \{T^{-1} (u) = 1_X\} \cup \{T^{-1}(u_i); i \in I\} \cup \{(T^{-1}(u_i))^2;i \in I\} \subset C(X).
$$
Der von dieser Teilmenge erzeugte Unterraum $G$ von $C(X)$  besitzt die Eigenschaft, dass
$$
\overline{\mathfrak{A}(Q)}^{C(X)} \subset \hat{G}_X \subset C(X) .
$$
Hierbei bezeichnet $\overline{\mathfrak{A} (Q)}^{C(X)}$ die kleinste abgeschlossene Teilalgebra, die $Q$ enth\"alt.
Nun ist $\{u_i; i \in I\}$ ein Erzeugendensystem von $E$, also auch von $E_u$. Diese Sprechweise ist sinnvoll, weil $E_u$ selbst ein lokalkonvexer Vektorverband ist. \\
$\{T^{-1} (u_i); i \in I\} \subset \overline{\mathfrak{A} (Q)}$ ist also ein Erzeugendensystem des lokalkonvexen Vektorverbandes $C(X)$. Damit  ist die abgeschlossene Algebra $\overline{\mathfrak{A} (Q)}$ ein abgeschlossener Vektorunterverband von $C(X)$, der ein Erzeugendensystem von $C(X)$ enth\"alt, also gilt
$$
C(X) \subset \overline{\mathfrak{A}(Q)} \subset \hat{G}_X \subset C(X)
$$
und damit
$$
\hat{G}_X = C(X) ,
$$
also
$$
T(Q) \subset E_u = T(C(X)) = T(\hat{G}_X) .
$$
Nach Voraussetzung gilt weiter
$$
E = \overline{E_u}^E = \overline{T(\hat{G}_X)}^E .
$$
Ist nun $(T_n)_{n \in \mathbb{N}}$ eine gleichstetige Folge positiver linearer Abbildungen und $P : E \to F$ ein stetiger Verbandshomomorphismus mit
$$
T_nv \to Pv \;\; \mbox{in} \;\; F \;\;  \mbox{f\"ur alle} \;\; v \in \{u\} \cup \{u_i,u_i^{(2)}; i \in I\} ,
$$
so bedeutet dies ja
$$
T_n \circ T (q) \to P \circ T(q) \;\; \mbox{in}\;\; F \;\; \mbox{f\"ur alle} \;\; q \in Q .
$$
$Q$ enth\"alt $1_X$, also folgt mit Theorem 1:
$$
T_n f \to Pf \;\; \mbox{in}\;\; F \;\; \mbox{f\"ur alle} \;\; f \in \overline{T(\hat{G}_{supp \; P \circ T})}^E .
$$
Wegen
$$
E = \overline{T(C(X))}^E \supset \overline{T(\hat{G}_{supp \; P \circ T})}^E \supset \overline{T(\hat{G}_X)}^E = E
$$
folgt also
$$
T_n f \to Pf \;\; \mbox{in}\;\; F \;\; \mbox{f\"ur alle} \;\; f \in E .
$$
Also gilt die Aussage (ii). 

\medskip

\noindent Die Behauptung (iii) ergibt sich nun sofort aus (ii). \hfill $\square$

\bigskip \bigskip

\noindent {\bf Bemerkung:} Endlich erzeugte, folgenvollst\"andige lokal-konvexe Vektorverb\"ande $E$ sind insbesondere Gegenstand der Untersuchungen von {\sc wolff} \cite{wolff1, wolff2}. Als besonders interessant erweist es sich hier, dass in dieser Situation $E$ automatisch einen quasiinneren Punkt $u$ enth\"alt. Ist etwa $A = \{u_1, \ldots , u_k\}$ ein Erzeugendensystem von $E$ (d.h. $E$ ist der kleinste abgeschlossene Vektorunterverband, der $A$ enth\"alt), so setze man $u : = \sum_{i=1}^k |u_i|$. Dann ist $u$ ein quasiinnerer Punkt des positiven Kegels $E^+$ von $E$, d.h. $E_u$ ist dicht in $E$. $E_u$ ist ja ein $A$ enthaltender Vektorunterverband von $E$. 

Es gelingt Wolff \cite{wolff3}, endlich erzeugte Banachverb\"ande vollst\"andig zu charakterisieren und als gewisse Funktionenr\"aume (sogenannte Banach-\-Funktionenverb\"ande) \"uber geeigneten kompakten Teilmengen des $\mathbb{R}^n, n \ge 1$, darzustellen.

Im Spezialfall sogenannter minimaler und separabler Banachverb\"ande gelingt es dabei sogar, auf ein Erzeugendensystem von nur zwei Elementen zu schlie\ss{}en.

\noindent
 $\begin{array}{ll}
\textrm{Heiner Gonska}\\
 \textrm{University of Duisburg-Essen} \\
 \textrm{Faculty of Mathematics} \\
 \textrm{D-47048 Duisburg} \\
 \textrm{Germany}\\
\textrm{e-mail: heiner.gonska@uni-due.de} \end{array}
$


\end{document}